\magnification\magstephalf
\baselineskip14pt
\parskip3pt

\def\proof{\noindent{\bf Proof. }}
\def\EC{{\rm EC}}
\def\OC{{\rm OC}}
\def\E{{\rm E}}
\def\O{{\rm O}}
\def\mod{\rm mod}
\def\pfbox
  {\hbox{\hskip 10pt\lower2pt\vbox{\hrule
  \hbox to 7pt{\vrule height 5pt\hfill\vrule}
  \hrule}}\hskip3pt}
\def\meno{\medskip\noindent}
\def\bib[#1] {\par\noindent\hangindent 40pt\hbox to20pt{[#1]\hfil}}

\def\Borch{1}
\def\Cul{2}
\def\CDGT{3}
\def\CDS{4}
\def\CG{5}
\def\EKLP{6}
\def\GM{7}
\def\Krew{8}
\def\Lov{9}
\def\MaMi{10}
\def\Stan{11}
\def\Wei{12}

\centerline{\bf Aztec diamonds, checkerboard graphs, and spanning trees}
\bigskip
\centerline{Donald E. Knuth, Stanford University}
\bigskip\bigskip
{\narrower\smallskip\noindent
{\bf Abstract.} 
This note derives the characteristic polynomial of a graph that represents
nonjump moves in a generalized game of checkers. The number of spanning trees
is also determined.
\smallskip}

\bigskip\noindent
Consider the graph on $mn$ vertices $\{(x,y)\mid 1 \leq x\leq m$, 
$1 \leq y\leq n\}$, with $(x,y)$ adjacent to $(x',y')$ if and only if
$\vert x-x'\vert=\vert y-y'\vert=1$. This graph consists of disjoint subgraphs
$$\eqalign{\EC_{m,n}&=\{(x,y)\mid x+y\hbox{ is even}\}\,,\cr
\OC_{m,n}&=\{(x,y)\mid x+y\hbox{ is odd}\}\,,\cr}$$
having respectively $\lceil mn/2\rceil$ and $\lfloor mn/2\rfloor$ vertices.
When $mn$ is even, $\EC_{m,n}$ and $\OC_{m,n}$ are isomorphic.
The special case $\OC_{2n+1,2n+1}$ has been called an {\it Aztec diamond
of order\/}~$n$ by Elkies, Kuperberg, Larsen, and Propp 
[\EKLP],
who gave several interesting proofs that it contains exactly $2^{n(n+1)/2}$
perfect matchings. Richard Stanley recently conjectured
[\Stan]
that $\OC_{2n+1,2n+1}$ contains exactly 4~times as many spanning trees as
$\EC_{2n+1,2n+1}$, and it was his conjecture that motivated the present note.
We will see that Stanley's conjecture follows from some even more remarkable
properties of these graphs.

In general, if $G$ and $H$ are arbitrary bipartite graphs having parts of
respective sizes $(p,q)$ and $(r,s)$, their {\it weak direct product\/}
$G\times H$ has $(p+q)(r+s)$ vertices $(u,v)$, with $(u,v)$ adjacent to
$(u',v')$ if and only if $u$ is adjacent to~$u'$ and $v$ to~$v'$. This graph
$G\times H$ divides naturally into even and odd subgraphs
$$\eqalign{\E(G,H)&=\{\,(u,v)\mid u\in G\hbox{ and }v\in H\hbox{ are in 
corresponding parts}\}\,,\cr
\O(G,H)&=\{\,(u,v)\mid u\in G\hbox{ and }v\in H\hbox{ are in opposite 
parts}\}\,,\cr}$$
which are disjoint. Notice that $\E(G,H)$ and $\O(G,H)$ have $pr+qs$ and $ps+qr$
vertices, respectively. Our graphs $\EC_{m,n}$ and $\OC_{m,n}$ are just
$\E(P_m,P_n)$ and $\O(P_m,P_n)$, where $P_n$ denotes a simple path on $n$~points.

Let $P(G;x)$ be the characteristic polynomial of the adjacency matrix of a
graph~$G$. The eigenvalues of $\E(G,H)$ and $\O(G,H)$ turn out have a simple
relation to the eigenvalues of~$G$ and~$H$:

\proclaim
Theorem 1. The characteristic polynomials\/ $P\bigl(\E(G,H);x\bigr)$ and
$P\bigl(\O(G,H);x\bigr)$ satisfy
$$P\bigl(\E(G,H);x\bigr)\,P\bigl(\O(G,H);x\bigr)
=\prod_{j=1}^{p+q}\;\prod_{k=1}^{r+s}\;(x-\mu_j\lambda_k)\,;\eqno(1)$$
\vskip-10pt
$$P\bigl(\E(G,H);x\bigr)=x^{(p-q)(r-s)}\,P\bigl(\O(G,H);x\bigr)\,.\eqno(2)$$

\proof
This theorem is a consequence of more general results proved in [\GM], as
remarked at the top of page 67 in that paper, but for our purposes a
direct proof is preferable.

Let $A$ and $B$ be the adjacency matrices of $G$ and $H$. It is well known 
[\Cul;~\Wei] that the adjacency matrix of $G\times H$ is the Kronecker product
$A\otimes B$, and that the eigenvalues of $A\otimes B$ are $\mu_j\lambda_k$
when $A$ and~$B$ are square matrices having eigenvalues~$\mu_j$
and~$\lambda_k$, 
respectively~[\MaMi, page~24]. Since the left side of~(1) is just
$P(G,H;x)$, equation~(1) is therefore clear.

Equation (2) is more surprising, because the graphs $\E(G,H)$ and $\O(G,H)$ often
look completely different from each other. But we can express $A$ and~$B$ in
the form 
$$A=\pmatrix{O_p&C\cr \noalign{\smallskip} C^T&O_q\cr}\,,\qquad
B=\pmatrix{O_r&D\cr \noalign{\smallskip} D^T&O_s\cr}\,,\eqno(3)$$
where $C$ and $D$ have respective sizes $p\times q$ and $r\times s$, and where
$O_k$ denotes a $k\times k$ matrix of zeroes. It follows that the adjacency
matrices of $\E(G,H)$ and $\O(G,H)$ are respectively
$$\pmatrix{O_{pr}&C\otimes D\cr \noalign{\smallskip}
C^T\otimes D^T&O_{qs}\cr}\qquad
{\rm and}\qquad
\pmatrix{O_{ps}&C\otimes D^T\cr \noalign{\smallskip}
C^T\otimes D&O_{qr}\cr}\,.\eqno(4)$$
We want to show that these matrices have the same eigenvalues, except for the
multiplicity of~0.

One way to complete the proof is to observe that the $k$\/th powers of both
matrices have the same trace, for all~$k$. When $k=2l$ is even, both matrix
powers have trace $\bigl({\rm tr}(CC^T)^l+{\rm tr}(C^TC)^l\bigr)
\bigl({\rm tr}(DD^T)^l+{\rm tr}(D^TD)^l\bigr)$ by
[\MaMi, pages 8, 18]; and when
$k$ is odd the traces are zero. The coefficients $a_1,a_2,\ldots$ of
$P(G;x)=x^{\vert G\vert}(1-a_1x^{-1}+a_2x^{-2}-\cdots\,)$ are completely
determined by the traces of powers of the adjacency matrix of any graph~$G$,
via Newton's identities; therefore (2) holds.\quad\pfbox

\proclaim
Corollary 1. The characteristic polynomials $P(\EC_{m,n};x)$ and
$P(\OC_{m,n};x)$ satisfy
$$P(\EC_{m,n};x)\,P(\OC_{m,n};x)=\prod_{j=1}^m\;\prod_{k=1}^n
\left(x-4\cos\,{j\pi\over m+1}\,\cos\,{k\pi\over n+1}\right)\,;\eqno(5)$$
\vskip-10pt
$$P(\EC_{m,n};x)=x^{mn\,\mod\, 2}\,P(\OC_{m,n};x)\,.\eqno(6)$$

\proof
It is well known 
[\Lov, problem 1.29; or \CDS, page 73],
that the eigenvalues of the path graph~$P_m$ are
$$\left\{2\cos\,{\pi\over m+1}\,,\;2\cos\,{2\pi\over m+1}\,,\,\ldots\,,
\;2\cos\,{m\pi\over m+1}\,\right\}\,.\eqno(7)$$
Therefore (1) and (2) reduce to (5) and (6).\quad\pfbox

\proclaim
Theorem 2. If $m\geq 2$ and $n\geq 2$, the number of spanning trees of\/
$\EC_{m,n}$ is $P(\OC_{m-2,n-2};4)$, and the number of spanning trees of\/
$\OC_{m,n}$ is\/ $P(\EC_{m-2,n-2};4)$.

\proof
Both $\EC_{m,n}$ and $\OC_{m,n}$ are connected planar graphs, so they have
exactly as many spanning trees as their duals
[\Lov, problem 5.23].
The dual graph $\EC_{m,n}^{\ast}$ has vertices $(x,y)$ where $1<x<m$ and
$1<y<n$ and $x+y$ is odd, corresponding to the face centered at $(x,y)$; it
also has an additional vertex~$\infty$ corresponding to the exterior face. All
its non-infinite vertices have degree~4, and when $\EC_{m,n}^{\ast}$ is
restricted to those vertices it is just $\OC_{m-2,n-2}$. Therefore the
submatrix of the Laplacian of $\EC_{m,n}^{\ast}$ that we obtain by omitting
row~$\infty$ and column~$\infty$ is just $4I-M$, where $M$ is the adjacency
matrix of $\OC_{m-2,n-2}$. And the number of spanning trees of
$\EC_{m,n}^{\ast}$ is just the determinant of this matrix, according to the
Matrix Tree Theorem
[\Borch; \Lov, problem 4.9; \CDS, page 38].

A similar argument enumerates the spanning trees of $\OC_{m,n}$. The basic idea
of this proof is due to Cvetkovi\'c and Gutman
[\CG];
see also
[\CDGT, pages 85--88]. \pfbox

\medskip
Combining Theorem 2 with equation (6) now yields a generalization of Stanley's
conjecture~[\Stan].

\proclaim
Corollary 2. When $m$ and $n$ are both odd, $\OC_{m,n}$ contains exactly 4
times as many spanning trees as $\EC_{m,n}$. \pfbox

Another corollary that does not appear to be obvious 
a~priori follows from Theorem~2 and equation~(5):

\proclaim
Corollary 3. When $m$ and $n$ are both even, $\EC_{m,n}$ contains an odd number
of spanning trees.

\proof
The adjacency matrix of $P_m$ is nonsingular $\mod\;2$ when $m$ is even.
Hence the adjacency matrix of  $\EC_{m,n}\cup\OC_{m,n}$ is
nonsingular $\mod~2$. Hence $P(\EC_{m,n};4)\equiv 1\;(\mod~2)$. \pfbox

Stanley
[\Stan]
tabulated the number of spanning trees in $\OC_{2n+1,2n+1}$ for $n\leq 6$ and
observed that the numbers consisted entirely of small prime factors. For
example, the Aztec diamond graph $\OC_{13,13}$ has exactly $2^{32}\cdot 
3^7\cdot
5^5\cdot 7^3\cdot 11^3\cdot 13^2\cdot 73^2\cdot 193^2$ spanning trees. One way
to account for this is to note that the number of spanning trees in
$\OC_{2n+1,2n+1}$ is
$$\eqalignno{&4^{2n-1}\,\prod_{j=1}^{n-1}\,\prod_{k=1}^{n-1}
\left(4-4\cos\,{j\pi\over 2n}\,\cos\,{k\pi\over 2n}\right)
\left(4+4\,\cos\,{j\pi\over 2n}\,\cos\,{k\pi\over 2n}\right)\cr
&\qquad =4^{2n-1}\,\prod_{j=1}^{n-1}\,\prod_{k=1}^{n-1}
\bigl(4-(\omega^j+\omega^{-j})(\omega^k+\omega^{-k})\bigr)
\bigl(4+(\omega^j+\omega^{-j})(\omega^k+\omega^{-k})\bigr)\,,&(8)\cr}$$
where $\omega=e^{\pi i/2n}$ is a primitive $4n$\/th root of unity.
Thus each factor $4-(\omega^j+\omega^{-j})(\omega^k+\omega^{-k})$ is an
algebraic integer in a cyclotomic number field, and all of its conjugates
$4-(\omega^{jt}+\omega^{-jt})(\omega^{kt}+\omega^{-kt})$ appear. Each product
of conjugate factors is therefore an integer factor of~(8).

Let us say that the edge from $(x,y)$ to $(x',y')$ in the graph is positive or
negative, according as $(x-x')(y-y')$ is $+1$ or~$-1$. The authors of
[\EKLP] showed that the generating function for perfect matchings in
$\OC_{2n+1,2n+1}$ is $(u^2+v^2)^{n(n+1)/2}$, in the sense that the coefficient
of~$u^kv^l$ in this function is the number of perfect matchings with
$k$~positive edges and $l$~negative ones. It is natural to consider the
analogous question for spanning trees: What is the generating function for
spanning trees of $\EC_{m,n}$ and $\OC_{m,n}$ that use a given number of
positive and negative edges? A~careful analysis of the proof of Theorem~2 shows
that the generating function for cotrees (the complements of spanning trees) in
$\OC_{m,n}$ is $P(\EC_{m-2,n-2};2u+2v)$, where $P$ now represents the
characteristic polynomial of the weighted adjacency matrix with positive and
negative edges represented respectively by~$u$ and~$v$. There are $\lceil
(m-1)(n-1)/2\rceil$ positive edges and ${\lfloor(m-1)(n-1)/2\rfloor}$ negative
edges altogether, so we get the generating function for trees instead
 of cotrees
by replacing~$u$ and~$v$ by~$u^{-1}$ and~$v^{-1}$, then multiplying by
$u^{\lceil(m-1)(n-1)/2\rceil}v^{\lfloor(m-1)(n-1)/2\rfloor}$. 
A~similar approach
works for~$\EC_{m,n}$.

Unfortunately, however, 
the polynomial $P$ does not seem to simplify nicely for
general~$u$ and~$v$, as it does when $u=v=1$. In the case $m=n=3$, the results
look reasonably encouraging because we have
$$\eqalign{P(\EC_{3,3};x)&=x^3\bigl(x^2-2(u^2+v^2)\bigr)\,,\cr
P(\OC_{3,3};x)&=(x+u+v)(x-u-v)(x+u-v)(x-u+v)\,.\cr}$$
But when $n$ increases to 5 we get
$$\eqalign{P(\EC_{3,5};x)&=x^4\bigl(x^2-2(u^2+uv+v^2)\bigr)
\bigl(x^2-2(u^2-uv-v^2)\bigr)\,,\cr
P(\OC_{3,5};x)&=x\bigl(x^2-(u^2+v^2)\bigr)
\bigl(x^4-3(u^2+v^2)x^2+2(u^2-v^2)\bigr)\,.\cr}$$
The quartic factor of $P(\OC_{3,5};x)$ cannot be decomposed into quadratics
having the general form 
$\bigl(x^2-(\alpha u^2+\beta uv+\gamma v^2)\bigr)\bigl(x^2-(\alpha'u^2+
\beta'uv+\gamma'v^2)\bigr)$, so it is unclear how to proceed. Some
simplification may be possible, because additional 
factors do appear when we set
$x=2u+2v$:
$$\eqalign{P(\EC_{3,5};2u+2v)
&=64(u+v)^4(u^2+3uv+v^2)(u^2+5uv+v^2)\cr
\noalign{\smallskip}
P(\OC_{3,5};2u+2v)
&=4(u+v)^3(3u^2+8uv+3v^2)(3u^2+14uv+3v^2)\cr
\noalign{\smallskip}
P(\EC_{5,5};2u+2v)
&=32(u+v)^5(3u^2+8uv+2v^2)(2u^2+8uv+3v^2)\cr
&\qquad (2u^4+24u^3v+53u^2v^2+24uv^3+2v^4)\cr
\noalign{\smallskip}
P(\OC_{5,5};2u+2v)
&=5(u+v)^4(u^2+4uv+v^2)(3u^2+8uv+3v^2)\cr
&\qquad
(15u^2+10uv+v^2)(u^2+10uv+15v^2)\,.\cr}$$
However, these factors are explained by the symmetries of $\EC_{m,n}$ and
$\OC_{m,n}$.

\meno
{\bf Acknowledgements.}
I thank an anonymous referee for suggesting the present form of Theorem~1,
which is considerably more general (and easier to prove) than the special case
I~had observed in the first draft of this note.
Noam Elkies gave me the insight about algebraic conjugates, when I was
studying the related problem of spanning trees in grids~[\Krew].
After writing this note, I learned from Richard Stanley that
the eigenvalues of the adjacency matrices of $\OC_{2n+1,2n+1}$
and $\EC_{2n+1,2n+1}$ were independently discovered by Tim Chow.

\vfill\eject

\centerline{\bf References}
\medskip
\bib
[\Borch]
C. W. Borchardt, ``Ueber eine der Interpolation entsprechende Darstellung der
Eliminations-Resultante,'' {\sl Journal f\"ur die reine und angewandte
Mathematik\/ \bf 57} (1860), 111--121.

\bib
[\Cul]
Karel {\v C}ulik, ``Zur Theorie der Graphen,'' {\sl {\v C}asopis pro
P{\v e}stov\'an{\'\i} Matematiky\/ \bf 83} (1958), 133--155.

\bib
[\CDGT]
Drago{\v s} M. Cvetkovi\'c, Michael Doob, 
Ivan Gutman, and Aleksandar Torga{\v s}ev,
{\sl Recent Results in the Theory of Graph Spectra}, Annals of Discrete
Mathematics {\bf 36} (1988).

\bib
[\CDS]
Drago{\v s} M. Cvetkovi\'c, Michael Doob, and Horst Sachs, {\sl Spectra of
Graphs\/} (New York: Academic Press, 1980).

\bib
[\CG]
D. Cvetkovi\'c and I. Gutman, ``A new spectral method for determining the
number of spanning trees,'' {\sl Publications de l'Institut Math\'ematique\/
\bf 29} (43) (Beograd, 1981), 49--52.

\bib
[\EKLP]
Noam Elkies, Greg Kuperberg, Michael Larsen, and James Propp,
``Alternating-sign matrices and domino tilings,''
{\sl Journal of Algebraic Combinatorics\/ \bf 1} (1992), 111--132 and 219--234.

\bib
[\GM]
C. Godsil and B. McKay, ``Products of graphs and their spectra,'' in
{\sl Combinatorial Mathematics IV}, edited by A. Dold and B. Eckmann,
{\sl Lecture Notes in Mathematics\/ \bf560} (1975), 61--72.

\bib
[\Krew]
Germain Kreweras, ``Complexit\'e et circuits eul\'eriens dans les sommes
tensorielles de graphes,'' {\sl Journal of Combinatorial Theory\/ \bf B24}
(1978), 202--212.

\bib
[\Lov]
L\'aszl\'o Lov\'asz, {\sl Combinatorial Problems and Exercises}, 2nd edition
(Budapest: Akad\'emiai Kiad\'o, 1993).

\bib
[\MaMi]
Marvin Marcus and Henrik Minc, {\sl A Survey of Matrix Theory and Matrix
Inequalities\/} (Boston: Allyn and Bacon, 1964).

\bib
[\Stan]
Richard P. Stanley, ``Spanning trees of Aztec diamonds,'' open problem
presented at a DIMACS meeting on Formal Power Series and Algebraic
Combinatorics, Piscataway, NJ, May 23--27, 1994.

\bib
[\Wei]
Paul M. Weichsel, ``The Kronecker product of graphs,'' {\sl Proceedings of the
American Mathematical Society\/ \bf 13} (1962), 47--52.

\bye